\documentclass[12pt]{article}
\usepackage{amsfonts}
\usepackage{amsmath}
\usepackage{amsthm}
\usepackage{amssymb}
\usepackage[english]{babel}
\usepackage{graphics}
\usepackage{latexsym}
\usepackage{makeidx}
\usepackage{longtable}
\usepackage{graphicx}

\usepackage{longtable}

\swapnumbers

\newcommand{\R}{\mathbb{R}} \newcommand{\C}{\mathbb{C}}
\newcommand{\N}{\mathbb{N}} 
\newcommand{\larrow}{\longrightarrow}  \newcommand{\lmapsto}{\longmapsto}
\newtheorem{theorem}{Theorem}[section]
\newtheorem{corollary}[theorem]{Corollary}
\newtheorem{lemma}[theorem]{Lemma}
\newtheorem{proposition}[theorem]{Proposition}

\theoremstyle{definition}
\newtheorem{definition}[theorem]{Definition}
\newtheorem{remark}[theorem]{Remark}
\newtheorem{remarks}[theorem]{Remarks}

\begin{document}

\pagestyle{myheadings} \markboth{\centerline{\small{\sc Nikolaos
I. Katzourakis}}}
         {\centerline{\small{\sc The Logarithmic Riemann Surface in a Geometrical Framework}}}

\title{The Riemann Surface of the Logarithm Constructed in a Geometrical
    Framework}

\author{ Nikolaos I. Katzourakis\footnote{Address: Troias 6, 18541, Kaminia, Pireaus,
Greece}\\[8pt]
         {University of Athens, Department of Mathematics}\\[8pt]
         { \small e-mail: nkatzourakis@math.uoa.gr}}

%\address{Troias 6, 18541, Kaminia, Pireaus, Greece}

%\keywords{Logarithmic Riemann Surface}

\date{}

\maketitle

%%% ----------------------------------------------------------------------

\begin{abstract}

\noindent The logarithmic Riemann surface $\Sigma_{log}$ is a
classical holomorphic $1$-manifold. It lives into $\R^4$ and
induces a covering space of $\C \smallsetminus \{0\}$ defined by
$\exp_{\C}$.

This paper suggests a geometric construction of it, derived as the
limit of a sequence of \emph{vector fields extending $\exp_{\C}$}
suitably to \emph{embeddings of $\C$ into $\R^3$}, which turn to
be \emph{helicoid} surfaces living into $\mathbb{C} \times
\mathbb{R}$. In the limit we obtain a bijective complex
exponential on the covering space in question, and thus a
well-defined complex logarithm. In addition, the \emph{helicoids}
are \emph{diffeomorphic (not bi-holomorphic) copies} of
$\Sigma_{log}$ as \emph{$C^{\infty}$-realizations living into
$\R^3$}, without obstruction.

Our approach is purely \emph{geometrical} and \emph{does not
employ any tools provided by the complex structure}, thus
\textit{holomorphy is no longer necessary to obtain
constructively this Riemann surface} $\Sigma_{log}$. Moreover, the
differential geometric framework we adopt affords explicit
generalization on submanifolds of $\C^m \times \R^m$ and certain
corollaries are derived.


\medskip \noindent
{\it MSC Subject Classification $2000$:} 33B30, 33B10, 37E30,
53A05, 53A07.
\end{abstract}

$\hspace{180pt}\underset{\scriptstyle{\emph{unless you have more
than one proof for it"\hspace{20pt}}}}{\scriptstyle{\emph{"You can
not claim you are well aware of a theorem}}}$

\medskip $\hspace{300pt}\scriptstyle{\textit{M. Atiyah's dictum}}$

\section*{Introduction.}

Riemann surfaces have been, at least conceptually, first
introduced by Riemann in his celebrated 1851 PhD dissertation. He
considered a surface spread over $\C$ with several sheets lying
over it, in such a way that a complex multi-valued function with
branches becomes a true well defined function, mapping each
branch to a sheet. His posterity gave a rigorous definition of
the (so called) \emph{Concrete Riemann surfaces} as covering
spaces of $\C$ (for an enlightening review see \cite{Remm}).

Klein was the first to put away the covering space approach and
adopted a differential geometric one: he studied complex functions
living on a curved surface of some ambient Euclidean space. In
fact, it was the first time the \emph{Atlas of holomorphic
structure} was introduced. At those times, "Riemann surface" meant
compact $2$-manifold with an arc-length element $ds^2$ and
bi-holomorphic transition functions in the atlas. Several years
later, adding Cantor's $2$nd countability and Hausdorff's axioms
(not before the 1920's), Rad\'o reaches via triangulations the
\emph{Abstract Riemann surfaces'} definition: a \emph{Hausdorff
$2$nd-countable topological surface with complex structure}.

The equivalence of \emph{Concrete} and \emph{Abstract} Riemann
surfaces follows from a theorem of Behnke \& Stein ($1947-49$)
\cite{Behnke-Stein}, improved by Gunning and Narasimhan (1967)
\cite{Gunning-Nar}.

In modern times, the research interest has been transferred to
the study of their \emph{Moduli (\& Teichm\"uller) spaces}, being
in effect sets of equivalence classes of distinct complex
structures, modulo the action of orientation preserving
diffeomorphism of the surface with certain topological structures
\cite{L-S-Y}, as well as to their vasty applications in
\emph{theoretical physics}, specifically \emph{$M$-theory
unifying the several String theories} (see e.g. \cite{Krasnov}).

In this paper we follow a semi-classical approach, standing on
the line between the \emph{abstract} and the \emph{concrete}. Our
principal result is that \emph{we recover in a pure
differential-geometric fashion results classically obtained via
the holomorphic structure itself}, \textit{bypassing analytic
continuation}.

Technically, the employed apparatus consists of a sequence of
vector fields $\{{\rm Exp}_n\}_{n \in \N}$ defined on $\C$ and
valued in (the tangent bundle of) $\R^3$. Each field ${\rm Exp}_n
: \; \C \larrow \C \smallsetminus \{0\} \times \R$ (\emph{Def.}
\ref{Exp definition}) geometrically is in fact a smooth parametric
surface of $\R^3$ called (exponential) helicoid (\emph{Def.}
\ref{Helicoid def}). These ${\rm Exp}_n$'s extend (due to the
extra component) the ordinary $\exp_{\C}$ in such a way that they
become $C^{\infty}$-diffeomorphisms onto their images (\&
embeddings of $\C$ into $\C \times \R$) removing $2 \pi i $
periodicity.

The ${\rm Exp}_n$-maps by construction converge to the covering
space of the punctured plane $\exp_{\C} : \; \C \larrow \C
\smallsetminus \{0\}$. Bijectivity is preserved in the limit
allowing a well defined limit logarithm $\log_{\C}$ to be
introduced on the direct limit surface $\underset{\larrow}{\lim}
\;{\rm Exp}_n (\C) \cong \exp_{\C}(\C)$ (in \ref{Convergence to
well-defined Log}) which essentially is its \emph{Riemann surface
$\Sigma_{log}$} immersed into the punctured plane. Even though all
the countably many sheets are "compressed" to a single sheet, we
still have the correspondence of the $\log_{\C}$-branches and the
sheets of the covering space via the convergence mechanism.

Notwithstanding, we show that each ${\rm Exp}_n$ is
$C^{\infty}$-diffeomorphic to $\Sigma_{log}$, when the latter is
considered as equipped with the smooth sub-atlas of the
holomorphic one coming from the ambient $\R^4$-space. The fact
that the helicoids can only be diffeomorphic copies of
$\Sigma_{log}$ and not bi-holomorphic is imposed by topological
obstructions due to the (general) non-imbeddability of Riemann
surfaces into $\R^3$.

The last \emph{Section} is devoted to high-dimensional
generalizations on multi-helicoid submanifolds of Euclidean
$3m$-space.

In view of the above, a natural question perhaps arises to the
reader: what kind of outstanding property do the helicoids enjoy
and what is so special about the {\rm Exp} vector fields? The
answer is that they are an \emph{ad hoc} choice in order to
provide the desired properties via convergence. There may very
well exist even more privileged surfaces of $\R^3$, but the smooth
realizations of the holomorphic $\Sigma_{log}$ into $3$-space are
\emph{dimension-wise} clearly \emph{optimal}.

\section{Preliminaries.}

We collect a few preparatory results which will be needed in the
main part of the paper, noted here for the reader's convenience.

\begin{definition} \textbf{(Helicoid)} \label{Helicoid def}
The smooth (parametric) surface of the \emph{(exponential)
helicoid} (see e.g. \cite{ON}, \cite{Sharipov}), is the map $X
:\; \R^2 \larrow \R^3$ given by
\[
    (u, v) \lmapsto (e^u \cos v,\; e^u \sin v, \;a v)\;,\;\;\; \  a >0
\]

where $a$ is a parameter. We consider this as the (global)
coordinate system of a $2$-manifold imbedded in $\R^3$. In
complex coordinates we may write
     $$X\;:\;\C \larrow \C \times \R \hspace{10pt}\;\;:\;\;$$
     $$ u+i v \;\lmapsto\;
     (\exp_{\C} (u + i v),\; a v)\; \equiv
     \; (\exp_{\C}(z),\; a {\rm Im}(z))$$

where $z \equiv u + i \; v$.The complex functions $Re$, $Im$,
$\exp_{\C}$ $\;:\;\C \larrow \C$ denote the "\emph{real part}",
the "\emph{complex part}" and the "\emph{complex exponential}" of
$\C$ respectively.
\end{definition}

The terminologies of imbeddings and immersions we follow are the
standard ones, referring e.g. to \cite{G-H-L}, \cite{K-N}.

We employ the convention that \emph{"imbedding"} stands for
topological imbedding and \emph{"embedding"} for geometrical
embedding in the sense of smooth manifolds.

In this paper \emph{smoothness} means \emph{$C^{\infty}$ -
smoothness} in the usual geometric sense and diffeomorphism
stands for $C^{\infty}$-diffeomorphism.

Furthermore, as usual by the term \emph{smooth $n$-manifold $M$}
(and in particular, surface) we mean a $C^\infty$-smooth
connected, paracompact Hausdorff manifold, of real $n$ dimensions
$\dim_{\R}(M) = n$ ($2$, in the 2nd case).

\begin{remark} \label{Euclidean vector fields} Let now
$$\exp_{\C}\; :\;\C \larrow \C \smallsetminus \{ 0\}\;\; :\;\;\; z
\lmapsto \exp_{\C} (z)\equiv e^z := \sum_{n=0}^{\infty}
\frac{z^n}{n !}$$

be the ordinary complex exponential function (\textit{the
notations $\exp_{\C}(z)$ and $e^z$ will be occasionally exchanged
without comments}). From the differential-geometric wiewpoint,
this map may be considered as a vector field, tangent on the
$2$-dimensional (flat) manifold $\C \cong \R^2$, i.e.
\[
    \C \equiv \R^2 \ni \; ue_1 +ve_2 \; \longmapsto \; (e^u \cos v)e_1 +
    (e^u \sin v)e_2 .
\]
\noindent This owes to the fact that if $\xi :\R^n \larrow \R^n$
is a vector field of $\R^n$ (in the \emph{vector calculus}
sense), where $x \lmapsto \xi (x)$, then the obvious
identification $\xi (x) \equiv (x, \xi (x))$ is the requested one
in order to consider $\xi$ valued in the tangent bundle $T\R^n
\cong \R^n \times \R^n \cong \R^{2n}$ (vector field in the
geometric sense) by global triviality of the vector (tangent in
our case) bundle (\cite{K-N}, \cite{G-H-L}, \cite{S}, etc).

The current task is to define a certain extension of the
exponential $\exp_{\C}$ from $\C$ to $\C \times \R$ that removes
the $2 \pi i$-periodicity and become injective due to extra real
component.
\end{remark}

These remarks introduce the idea of the following definition:

\begin{definition} \textbf{(The exponential field)} \label{Exp definition}
The exponential vector field on $\C$ is the map
\[
   {\rm Exp}_a \;\;:\;\;\; \C\; \larrow \; \C \times \R
\]
given by:
\[
  u + i v \; \lmapsto {\rm Exp}_a (u +i v)\;:=\;(\exp_{\C} (u +i v),\;a v)
\]
where $a > 0$ is a parameter to be fixed as will in the sequel.
Operationally the map can given as
\[
  {\rm Exp}_a \;=\;( \exp_{\C} , \; a \; {\rm Im})
\]
\end{definition}

The following technical fact shows that this vector-wise
extension of the $\exp_{\C}$ function on $\C$ is the requested
one that provides injectivity of the map, now considered as a
vector field.

\begin{proposition} \textbf{(Structural Properties of {\rm Exp})} \label{Structural Properties of Exp}

\noindent The map $\;{\rm Exp}_a :\; \C \larrow \C \times \R$
given by \ref{Exp definition} has the following properties:

\emph{a)} It is a vector field defined on the (trivially) embedded
submanifold $\C \equiv \C \times \{0\} \hookrightarrow \C \times
\R \cong \R^3$ and valued into the (tangent bundle of the)
ambient space $\R^3$.

\emph{b)} It is a diffeomorphism onto its image, thus an embedding
of the plane into $\C \times \R$ and, in particular, invertible.
\end{proposition}

\noindent As our context suggests, we consider submanifolds as
subsets of $\R^n$, where the inclusion map is an imbedding and
the Atlas that determines differential structure and topology is
the induced from $\R^n$.

\begin{proof}
\emph{a)} We recall that a vector field defined on a
\emph{submanifold} $M$ of $\R^n$ in not necessarily tangent to
the manifold in the usual differential geometric sense, but the
splitting of tangent spaces
            $$\R^n \; \cong \; T_{p}\R^n \; \cong
             \; T_{p}M \bigoplus T_{p}^{\perp}M
            \;,\;\;\;\;p \in \R^n$$

implies that it can be normal, or generally valued in $T \R^n
\cong \R^{2n}$. The very definition of ${\rm Exp}_a$ implies that
it maps a point $z$ of $\C$ to a vector ${\rm Exp}_a (z)$ of $\C
\times \R$
\[
  \C \; \ni \;\; z \;\longmapsto \;{\rm Exp}_a (z) \;\;\in \; \C \times \R
\]
and \emph{Remark} \ref{Euclidean vector fields} implies that the
map is a vector field of $\C$ valued in the tangent bundle of $\C
\times \R$.

\emph{b)} If ${\rm Exp}_a(u_1 +i v_1,\; a v_1) = {\rm Exp}_a(u_2
+i v_2,\; a v_2)$ for $u_1, u_2, v_1, v_2 \in \R$, then this
amounts to
\[
                 (\exp_{\C}(u_1 +i v_1),\;a v_1)
              = (\exp_{\C} (u_2 +i v_2),\;a v_2)
\]
that is
\[
    \exp_{\R}(u_1) \exp_{\C}(i v_1) \;=\;\exp_{\R}(u_2)
     \exp_{\C}(i v_2) \;,\;\;\;\;\;v_1 = v_2
\]
and consequently $u_1 = u_2$ by the monotonicity of the real
exponential. This shows the injectivity of {\rm Exp}, thus it is
bijective onto its image ${\rm Exp}_a(\C)$.
$C^{\infty}$-differentiability of this map goes without saying,
since it has smooth component functions in $\C \times \R$.
Smoothness of the inverse map is also obvious (for the explicit
expression of the inverse map see \emph{Lemma} \ref{Expression of
Log} of the oncoming \emph{Section 2}). Consequently, the map is a
smooth (geometric) embedding of $\C$ into $\C \times \R$.
\end{proof}

\noindent Consider now the standard Euclidean (Riemannian) metric
$\delta$ of $\R^3$ with components the Kronecker deltas'
$\delta_{ab}$ (where $\delta_{aa} = 1$ and $0$ otherwise, $1 \leq
a \leq 3 $) (for the standard concepts of \emph{Riemannian
geometry} employed in this paper, we refer to \cite{DC},
\cite{G-H-L}, \cite{K-N}, \cite{ON}). Then, the (trivially)
embedded surface $\C$ gets a natural induced \emph{flat}
Riemannian metric, the standard inner product $< \cdot ,\cdot> \;
\equiv \; \delta$ of $\R^3$ by restriction of $\delta$ on $\C$,
which can be seen as the pull-back metric via the inclusion map:

    $$(\C,\; \delta |_{\C})\; \equiv \; (\C,\; i^{*}( \delta ))$$
    $$i\;:\;\;\C\; \hookrightarrow \;(\C \times \R, \; \delta )$$

In view of theses remarks, we obtain the next result giving some
geometric properties of the extended object ${\rm Exp}_a$ as a
vector field on the $(\C,\; i^{*}(\delta ))$.

\medskip \noindent \textbf{Notation.} $angl(\xi, \eta)(p)$ denotes
the angle of $2$ vector fields $\xi$, $\eta$ measured at the
point $p$ of the surface (manifold), with respect to its
Riemannian metric.

\begin{proposition} \textbf{(Geometric Properties of {\rm Exp})} \label{Geometric Properties of Exp}

\emph{a)} The vector field ${\rm Exp}_a : \C \larrow \C \times \R$
geometrically is a smooth exponential helicoid surface of the form
\ref{Helicoid def} and the following diagram commutes:

\[
\begin{picture}(300,100)
 \put(50,82){\vector(1,0){120}}
  \put(172,14){\vector(-1,0){105}}
%\put(50,78){\vector(3,-1){122}}
  \put(39,78){\makebox{$\C$}}
  \put(25,10){\makebox{$\C \smallsetminus \{0\}$}}
 \put(175,10){\makebox{$\C \smallsetminus \{0\} \times \R$}}
 \put(197,74){\vector(0,-1){56}}
  \put(43,76){\vector(0,-1){58}}
  \put(173,78){\makebox{${\rm Exp}_a(\C)$}}
\put(200,50){\makebox{$\scriptstyle{i}$}}
\put(47,50){\makebox{$\scriptstyle{\exp_{\C}}$}}
\put(110,85){\makebox{$\scriptstyle{{\rm Exp}_a}$}}
\put(109,18){\makebox{$\scriptstyle{\pi_{\C \smallsetminus
\{0\}}}$}}
 %\put(100,52){\makebox{$\scriptstyle{\exp_{\C}}$}}
\end{picture}
\]

that is, ${\rm Exp}_a$ commutes with the projection composed with
the inclusion of the surface $i : {\rm Exp}_a (\C)
\hookrightarrow \C \smallsetminus \{0\} \times \R$ and the
projection $\pi_{\C \smallsetminus \{0\}}\; :\;\; \C
\smallsetminus \{0\} \times \R \larrow \C \smallsetminus \{0\}$
\[
    (\;\pi_{\C \smallsetminus \{0\}} \circ \; i \;) \; \circ \;
     {\rm Exp}_a \;=\; \exp_{\C}
\]

\emph{b)} If the tangent plane $T_{u+i v} \C \cong \C \equiv \R^2$
is spanned by $e_1 = (1,0)$ and $e_2 = (0,1)$ and
               $$\xi_{AB} \equiv A e_1 + B e_2$$
is a general $2$-parameter family of tangent vectors, $A,\;B \in
\R$, $|A| + |B| \neq 0$, then

\[
   \cos [(\Phi_{AB})(u,v)] = \frac{e^{u}\;[A \cos v
   + B \sin v]} {\sqrt{e^{2u} + a^2 v^2} \;\; \sqrt{A^2 + B^2}}
\]

\medskip

where $\Phi_{AB}(u,v)$ is the angle of {\rm Exp} with $\xi_{AB}$
at $u+ iv$ in $\R^3$:
\[
(\Phi_{AB})(u, v) \equiv angl({\rm Exp}_a ,\;\xi_{AB})(u +i v)
\]
\end{proposition}

\begin{remark} The calculation of the aforementioned (non-constant) angle
shows that the ${\rm Exp}_a$-vector field is neither tangent nor
normal to $\C$.
\end{remark}

\begin{proof}
\emph{a)} The first claim goes without saying, just by comparing
the very definitions \ref{Helicoid def} and \ref{Exp definition}.
Thus, the map ${\rm Exp}_a$ is a vector field of $\R^3$ defined
on $\C$ and simultaneously its image constitutes a smoothly
embedded surface of $\R^3$, an exponential helicoid diffeomorphic
to $\C$.

Commutativity is a simple consequence of the form of the ${\rm
Exp}_a$-field:
\[
   (\pi_{\C \smallsetminus \{0\}} \circ \; i \;) \; \circ \;
     {\rm Exp}_a (z) \;=\; (\pi_{\C \smallsetminus \{0\}}\;
     \circ\; (\exp_{\C},\; a \; Im)) (z) \;=\; \exp_{\C}(z)
\]

\emph{b)} We recall the familiar Euclidean formula giving the
angle of $2$ vectors in $(\R^3, \delta)=(\R^3,\;<\cdot , \cdot>)$
\[
   \cos \; [ angl({\rm Exp}_a, \xi_{AB})(u, v)] \;=\; (\frac{<{\rm
   Exp}_a, \xi_{AB}>}{||{\rm Exp}_a||\; ||\xi_{AB}||}) (u + i v)
\]

The $2$-parameter vector (describing a general tangent vector) is
given due to the inclusion $\C \times \{0\} \hookrightarrow \C
\times \R$ by the formula
\[
  \xi_{AB} \;=\; A e_1 + B e_2 + 0 e_3 \;=\; A\; (1,0,0) + B\; (0,1,0) +
  (0,0,0) \;=\; (A,B,0)
\]
and provided that ${\rm Exp}_a = (e^u \cos v, e^u \sin v, a v)$ we
have
\[
   \cos [(\Phi_{AB})(u,v)] = \frac{A\; e^{u} \cos v
   + B\;e^u \sin v} {\sqrt{e^{2u} + a^2 v^2} \;\; \sqrt{A^2 + B^2}}
\]
which gives the requested formula and this completes the proof.
\end{proof}

\section{The Logarithmic Riemann surface.}

\textbf{(Part I) \underline{The Convergence Constructions}.}

\medskip \noindent The \emph{infinite-sheeted Riemann surface $\Sigma_{log}$
of the logarithm} $\log_{\C}$ (see \cite{B-N}, \cite{F-K},
\cite{Spr}), considered as a holomorphic manifold of (complex)
dimension $\dim_{ \C} (\Sigma_{log}) = 1$, lives into $\C^2 \cong
\R^4 $ and can be given as the graph of the function $e^z = w$
\[
                 \Sigma_{log}\; :=\; \{z \in \C \;/\; e^z = w\} \;
                 \subseteq \; \R^4
\]
or by the parametric equations
\[
                   \C \ni \; z \;\lmapsto \;(z , \exp_{\C}(z))\;
                    \in \C \times \C\; \cong\; \C^2
\]

$\Sigma^{\mathcal{O}}_{log}$ denotes the holomorphic manifold
$(\Sigma_{log}, \mathcal{A}_{\Sigma_{log}}^{\mathcal{O}})
\hookrightarrow (\R^4, id)$ with the induced analytic atlas
coming from $\R^4$. Classically, it is obtained by analytic
continuation, extending holomorphically $\exp_{\C}$ to "small"
discs along a disk $\mathbb{D}(0,1)$ of $\C$. Since the last
extension that overlaps with the first does not coincide with it,
we take a copy of $\C$ and proceed continuation to a next sheet
lying over the initial one. Continuing ad infinitum, we obtain a
covering space of $\C \smallsetminus \{0\}$ bi-holomorphic to
this surface.

Projecting to the second factor, we obtain the covering space map
on the punctured plane $\C \smallsetminus \{0\}$ (equivalence of
\emph{concrete} \& \emph{abstract} approach, \cite{Behnke-Stein},
\cite{Gunning-Nar}):
       $$\pi\;:\;\;\Sigma^{\mathcal{O}}_{log} \subseteq \C^2 \larrow \C
     \smallsetminus \{0\}\;\; : \;\;(z , \exp_{\C}(z)) \lmapsto \exp_{\C}(z)$$

In general there is no way to holomorphically imbed this
2-manifold into $\R^3$, unless we allow self-intersections, which
is no longer a homeomorphism onto its image but may still be an
immersion.

Our construction allows, as we shall see, to introduce a
well-defined complex logarithm. Note that \emph{all this apparatus
manages to bypass the holomorphic structure of this complex
manifold and no analytic continuation is anywhere used}.

\medskip

\noindent In this section we give the \emph{geometric}
construction of $\Sigma_{log}$. As stated in the
\emph{Introduction}, it is obtained as a special covering
manifold of $\C \smallsetminus \{ 0\}$ in the limit of a sequence
obtained by the exponential ${\rm Exp}_a$ images of $\C$.

Hence, we recall the map ${\rm Exp}_a :\; \C \larrow \C \times
\R$ given by \ref{Exp definition} and substitute the positive
parameter by the sequence $a_n = 1 / n$, $n \in \N$. Thus, we
obtain a sequence of vector fields (and surfaces of $\R^3$ as
well) $\{{\rm Exp}_n\}_{n \in \N}$, where
\[
    {\rm Exp}_n \;:\;\;\C \larrow \C \times \R \;\;:\;\;\; z \longmapsto \bigl(
    \exp_{\C}(z),\frac{1}{n} Im(z))\bigr) \;,\;\;\;\;\;\;n \in \N
\]

The following result is a well-known fact and can be found an any
elementary textbook, e.g. \cite{Mun}. For the definition of
\emph{covering manifolds} we refer to \cite{G-H-L}, \cite{S-T}.

\begin{lemma} \label{Covering space}
The pair $(\exp_{\C},\C)$, constitutes a covering manifold of $\C
\smallsetminus \{ 0\} $ with countably infinite covering sheets.
\end{lemma}

\medskip \noindent We are now in position to prove the main result
of this section. The sequence of maps $\{{\rm Exp}_n : \C \larrow
\C \times \R \}_{n \in \N}$ is a diffeomorphism onto its image, by
\textit{Proposition} \ref{Structural Properties of Exp}, for $a =
1/n$, $n\in \N$ and consequently a bijection.

We shall presently see that this property is preserved in the
limit $n \larrow \infty$. The limit surface is the covering space
produced by $\exp_{\C}\;:\; \C \larrow \C \smallsetminus \{ 0\}$
(\textit{Lemma} \ref{Covering space}), and introducing the
natural inclusion $\C \smallsetminus \{ 0\} \hookrightarrow \C
\smallsetminus \{ 0\} \times \R$ we may identify $\underset{n
\rightarrow \infty}{\lim}({\rm Exp}_n)$ on $\C$ with the usual
complex exponential, $\exp_{\C}$. Thus, we construct an
exponential injective on $\C$.

This situation resembles and simultaneously extends the analogous
case in elementary complex analysis (see e.g. \cite{B-N}), where
$\exp_{\C}$ is injective and invertible on the strips of $2\pi
i$-width only.

\begin{theorem} (\textbf{Convergence to Injective Complex Exponential})
 \label{The Riemann surface}

\noindent The sequence of diffeomorphisms $\{{\rm Exp}_n : \C
\overset{\cong}{\larrow} ({\rm Exp}_n)(\C) \subseteq \C
\smallsetminus \{0\} \times \R \}_{n\in \; \N}$ converges
(point-wise) to a diffeomorphism $\underset{n \rightarrow
\infty}{\lim} {\rm Exp}_n$ on the covering surface of $\C
\smallsetminus \{0\}$ produced by $\exp_{\C}\;:\;\C \;\larrow\;\C
\smallsetminus \{0\}$:
\[
    {\rm Exp}_n \; \underset{pw}{\xrightarrow{\;\;\;\;n \larrow
    \infty\;\;\;\;}}\;
    \exp_{\C} \times \{0\}\; \equiv \; \exp_{\C}
\]
and the lower triangle of the following diagram defines (in the
limit) a bijective complex exponential $\exp_{\C}$.
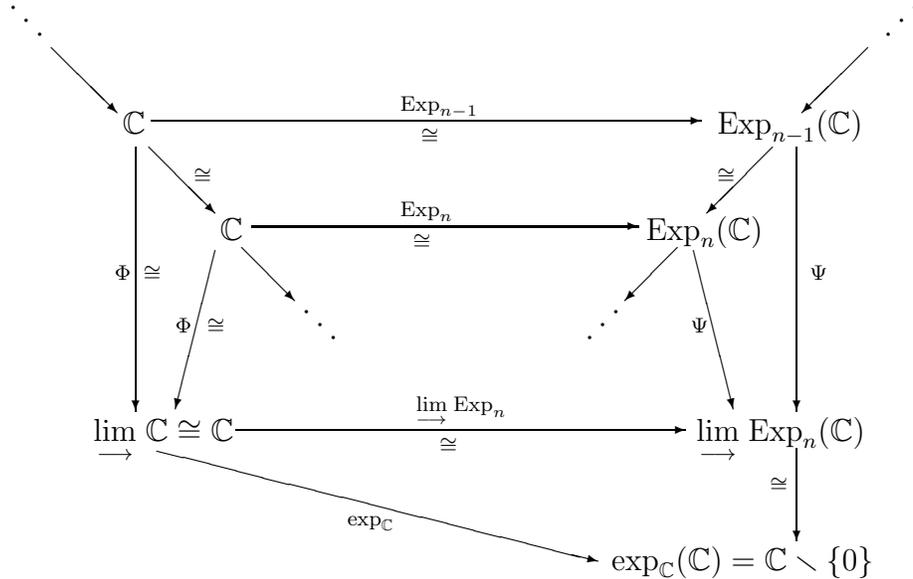
\[
\begin{picture}(300,220)
\put(250,165){\makebox{${\rm Exp}_{n-1}(\C)$}}
  \put(25,165){\makebox{$\C$}}
  \put(36,170){\vector(1,0){208}}

\put(250,147){\makebox{$\scriptstyle{\cong}$}}
 \put(52,147){\makebox{$\scriptstyle{\cong}$}}

\put(223,125){\makebox{${\rm Exp}_n(\C)$}}
  \put(62,125){\makebox{$\C$}}
  \put(74,130){\vector(1,0){146}}

 \put(235,122){\vector(-1,-1){20}}
\put(70,122){\vector(1,-1){20}} \put(93,95){\makebox{$\cdot$}}
\put(98,90){\makebox{$\cdot$}} \put(103,85){\makebox{$\cdot$}}

\put(210,95){\makebox{$\cdot$}} \put(205,90){\makebox{$\cdot$}}
\put(200,85){\makebox{$\cdot$}}

 \put(280,160){\vector(0,-1){100}}
  \put(30,160){\vector(0,-1){100}}
\put(280,46){\vector(0,-1){35}}
  %\put(30,46){\vector(0,-1){30}}

%\put(35,30){\makebox{$\scriptstyle{\cong}$}}
\put(270,30){\makebox{$\scriptstyle{\cong}$}}

 \put(271,160){\vector(-1,-1){25}}
\put(35,160){\vector(1,-1){25}}

\put(-8,200){\makebox{$\cdot$}} \put(313,200){\makebox{$\cdot$}}
\put(-13,205){\makebox{$\cdot$}} \put(318,205){\makebox{$\cdot$}}
\put(-18,210){\makebox{$\cdot$}} \put(323,210){\makebox{$\cdot$}}
 \put(-2,198){\vector(1,-1){25}}
 \put(307,198){\vector(-1,-1){25}}

 \put(60,121){\vector(-1,-4){15}}
 \put(241,121){\vector(1,-4){15}}

\put(22,110){\makebox{$\scriptstyle{\Phi}$}}
 \put(45,90){\makebox{$\scriptstyle{\Phi}$}}
\put(33,110){\makebox{$\scriptstyle{\cong}$}}
 \put(57,90){\makebox{$\scriptstyle{\cong}$}}
 \put(240,90){\makebox{$\scriptstyle{\Psi}$}}
\put(285,110){\makebox{$\scriptstyle{\Psi}$}}

%\put(110,9){\vector(1,0){128}}
 %\put(18,5){\makebox{$P(\C^{\#}) \cong \C \smallsetminus \{0\}$}}

 \put(210,0){\makebox{$\exp_{\C}(\C) = \C \smallsetminus \{0\}$}}

\put(68,53){\vector(1,0){170}}
 \put(14,49){\makebox{$\underset{\larrow}{\lim}\; \C \cong \C$}}
\put(242,49){\makebox{$\underset{\larrow}{\lim} \; {\rm
Exp}_n(\C)$}}
   \put(38,45){\vector(4,-1){167}}
\put(136,60){\makebox{$\scriptstyle{\underset{\larrow}{\lim}\;
{\rm Exp}_n}$}}
 \put(145,45){\makebox{$\scriptstyle{\cong}$}}

\put(130,134){\makebox{$\scriptstyle{{\rm Exp}_n}$}}
\put(130,174){\makebox{$\scriptstyle{{\rm Exp}_{n-1}}$}}
\put(135,122){\makebox{$\scriptstyle{\cong}$}}
\put(138,162){\makebox{$\scriptstyle{\cong}$}}
\put(110,17){\makebox{$\scriptstyle{\exp_{\C}}$}}
\end{picture}
\]
\end{theorem}

\begin{remarks}
\noindent Before giving the proof of this result which
essentially contains the requested construction, some explanatory
comments are necessary. First of all, the direct (inductive)
limits "$\underset{\larrow}{\lim}$" (for their definition in a
general \emph{Category} we refer to the classical textbook
\cite{Mac Lane}) are justified as follows:

The constant sequence of smooth manifolds
\[
  \C \; \underset{\cong}{{\xrightarrow{\;\;\;\; id \;\;\;\;\;}}}\; \C
  \; \underset{\cong}{{\xrightarrow{\;\;\;\; id \;\;\;\;\;}}} \;
  \C \; \larrow \; \cdots
\]

defines in a completely trivial way a direct limit set which is
$\C$, since $id \circ id = id $. In addition,
\[
  \underset{\larrow}{\lim} \; \C \; =\; \underset{n \rightarrow \infty}{\lim}
  \; \C \;=\; \C
\]
that is, the direct limit coincides with the usual sequence
convergence. Similarly, we have the sequence
\[
  {\rm Exp}_n (\C) \; \underset{\cong}{{\xrightarrow{\;\;\;\; \Theta_{n, n+1}
   \;\;\;\;\;}}}\; {\rm Exp}_{n+1} (\C)
  \; \underset{\cong}{{\xrightarrow{\;\;\;\; \Theta_{n+1, n+2} \;\;\;\;\;}}} \;
  {\rm Exp}_{n+2} (\C) \; \larrow \; \cdots
\]

where the diffeomorphisms $\Theta_{n,m}$, $n,\;m \in \N$ are
defined as
\[
  \Theta_{n,m} \; : \;\; {\rm Exp}_n (\C) \; \larrow \; {\rm Exp}_m (\C)
\]
by the formula
\[
  (\exp_{\C}(z),\; \frac{1}{n} Im (z)) \; \longmapsto \;
   \Theta_{n,m}\bigl(\exp_{\C}(z),\; \frac{1}{n} Im (z) \bigr) \;
\]
\[
  \hspace{140pt} :=\; (\exp_{\C}(z),\; [\frac{n}{m}] \frac{1}{n} Im (z))
\]

The direct limit properties can be easily verified, since by
construction of the $\Theta$-maps we get
\[
   \Theta_{n,m} \; \circ \; \Theta_{m,k} \;=\; \Theta_{m,k}
\]
as well as
\[
\Theta_{n,n} = id_{{\rm Exp}_n (\C)}
\]

This implies that
\[
  \underset{\larrow}{\lim} \; {\rm Exp}_n (\C) \; =\;
   \underset{n \rightarrow \infty}{\lim}\; {\rm Exp}_n (\C) \;=\;
   \exp_{\C}(\C) \times \{0\} \;\cong\; \C \smallsetminus \{0\}
\]

the $\Phi$-maps pictured in the diagram are the identity (up to
diffeomorphism) and the $\Psi$-maps are the projections on $\C
\smallsetminus \{0\}$ composed with the inclusions of ${\rm
Exp}_n (\C)$ into $\C \times \R$:
\[
  \Psi \;:=\; \pi_{\C \smallsetminus \{0\}} \circ i
\]
\end{remarks}

\begin{proof}
\textit{Proposition} \ref{Structural Properties of Exp} implies
that ${\rm Exp}_n$ is injective, $n \in \N$. Consequently, we have
\[
    {\rm Exp}_n (z) = {\rm Exp}_n (w) \;\; \Longrightarrow \;\; z = w
\]
Our task is to prove the following condition:
\[
    \exp_{\C} (z) \;=\; \exp_{\C} (w) \;\; \Longrightarrow \;\; z = w
\]

when $z,w$ live into $\C$ and by $\emph{a)}$ of \emph{Proposition}
\ref{Geometric Properties of Exp} we obtain
\[
{\rm Exp}_n \;{\xrightarrow{\;\;n \larrow \infty\;\;\;}}\;
\exp_{\C} \times \{0\} \equiv \exp_{\C}
\]
having used that $(1/n) \;{\xrightarrow{\;\;n \larrow
\infty\;\;\;}}\; 0$.

\noindent The well-known triangle inequality gives
\[
     \big\| {\rm Exp}_n (z) -{\rm Exp}_n
    (w) \big\|  \; \le  \;   \big\| {\rm Exp}_n (z) - (\exp_{\C}(z), 0) \big\| \;+
\]
\[
    +\; \big\|  {\rm Exp}_n (w) - (\exp_{\C}(w), 0) \big\| \;+\; | \exp_{\C}(w) -  \exp_{\C}(z)|
\]
When $n \larrow \infty$, we have
\[
    \big\| {\rm Exp}_n (z) - (\exp_{\C}(z), 0)
    \big\| \; {\xrightarrow{\;\;n \larrow \infty\;\;\;}} \; 0
\]
as well as
\[
    \big\|  {\rm Exp}_n (w) - (\exp_{\C}(w), 0)
    \big\| \; {\xrightarrow{\;\;n \larrow \infty\;\;\;}} \; 0
\]

\[
    \big\| {\rm Exp}_n (z) - {\rm Exp}_n (w) \big\| \; \leq \;
    |  \exp_{\C} (z)- \exp_{\C} (w) | \;
     {\xrightarrow{\;\;n \larrow \infty\;\;\;}} \; 0
\]
which means that in the limit $z = w$ by the injectivity assured
by \ref{Structural Properties of Exp}. This completes the proof.
\end{proof}

\begin{corollary}
\noindent The convergence is \underline{uniform} on the bounded
strips of $\C$ with bounded imaginary part $B\;\equiv \; \{z \in
\C \;/\; Im (z) < M,\; M > 0 \}$, that is, for finite spirals of
the ${\rm Exp}_n$-helicoids:
\[
    {\rm Exp}_n|_{B} \; \underset{U}{\xrightarrow{\;\;\;\;n
     \larrow \infty \;\;\;\;}}\; \exp_{\C} \times \{0\}|_{B} \;
      \equiv \; \exp_{\C}|_{B}
\]
for any $M > 0$.
\end{corollary}

\begin{remark}
The bounded strips of the form of $B$ corresponds in the limit to
a finite-sheeted covering space of $\C \smallsetminus \{0\}$:
\[
  \exp_{\C}|_{B} \;:\;\;B \; \larrow \; \C \smallsetminus \{0\}
\]
produced by the complex exponential.
\end{remark}

\begin{proof}
In order to see when the convergence of ${\rm Exp}_n$ to
$\exp_{\C}$ is uniform, let $\varepsilon > 0$. Then, \emph{Th.}
\ref{The Riemann surface} implies (using \emph{a)} of
\emph{Prop.} \ref{Geometric Properties of Exp})
\[
  \big\| {\rm Exp}_n (z) - (\exp_{\C}(z),\; 0) \big\| \;
   =  \; \big\| \bigl( \exp_{\C}(z), \;\frac{1}{n} Im (z)\bigr) - (\exp_{\C}(z),\; 0)
   \big\|  =\; \frac{1}{n} |Im (z) |
\]

\medskip \noindent Imposing $|Im (z)| < M$ bounded, we obtain
uniform convergence as claimed.
\end{proof}

\noindent The injectivity of the {\rm Exp} provided by
\textit{Proposition} \ref{Structural Properties of Exp} gives the
ability to introduce the inverse map (which we shall denote by
${\rm Log}$) as a well-defined generalization of the complex
logarithm:
\[
    {\rm Log}_n  \;\;:\;\;\; {\rm Exp}_n (\C)
    \; \subseteq \C \smallsetminus \{0\} \times \R \; \larrow \; \C
\]
where
\[
    {\rm Log}_n \;:=\; {\rm Exp}_n^{-1}\;,\;\;\;\;\;n \in \N
\]

\begin{lemma} \label{Expression of Log}
The inverse ${\rm Log}_n$ maps are given by
\[
   (K, L) \; \longmapsto \; {\rm Log}_{n}(K, \; L) \;=\; \ln(|K|)
   \;+\; i n L
\]
( $\ln$ denotes the real Napierian logarithm ).
\end{lemma}

\begin{proof}
The proof is a calculation. Set $$\exp_{\C}(z)\; \equiv \; K
\;\;,\;\;\;\;\;\; \frac{1}{n} Im (z) \; \equiv \; L$$ Then,
\[
  z \;=\; \log_{\C}(K) = \ln(|K|) + i \arg (K)
\]
Comparing the last relation with the formula $\log_{\C}(K) =
\ln(|Re (z)|) + i n L$ we conclude to the requested relation.
\end{proof}

We are now in position to introduce a well-defined complex
logarithm on the limit covering space produced by the exponential.

\medskip \noindent \textbf{Notational convention:} In the sequel we
adopt the following convention: $\Sigma_{log}^{C^{\infty}}$
denotes $\Sigma_{log}$ equipped with the induced
\underline{smooth} sub-atlas of the holomorphic one coming from
$\R^4$:
\[
   \Sigma_{log}^{C^{\infty}} \; \equiv \; (\Sigma_{log},\;
    \mathcal{A}_{\Sigma_{log}}^{C^{\infty}}) \; \hookrightarrow \; (\Sigma_{log},\;
    \mathcal{A}_{\Sigma_{log}}^{\mathcal{O}}) \; \subseteq \; (\R^4, Id)
\]

\begin{corollary} \label{Convergence to well-defined Log} \textbf{(Convergence
 to a well-defined $\log_{\C}$ on its Riemann Surface)}

The sequence of maps $({\rm Log}_n)_{n \in \N}$ of \ref{Expression
of Log} converges point-wise to a well-defined complex logarithm
and uniformly on the subsets of ${\rm Exp}_n (\C) \subseteq \C
\smallsetminus \{0\} \times \R$ of the form ${\rm Exp}_n (\C)
\bigcap \C \smallsetminus \{0\} \times (a,b)$:
\[
     {\rm Log}_n \; {\xrightarrow{\;\;\;n \larrow
     \infty\;\;\;\;}}\; (\exp_{\C} \times \{0\})^{-1}\;
      \equiv \; \log_{\C} |_{ \exp_{\C}(\C)}
\]

and the limit surface is diffeomorphic to the (immersed into $\C
\smallsetminus \{0\}$) \underline{Riemann} \underline{Surface} of
the \underline{logarithm}:
\[
    \hspace{20pt} \pi_{\C}(\Sigma_{log}^{C^{\infty}})\; \cong \; \lim_{\larrow}
     \;( {\rm Exp}_n (\C))\;=\; \exp_{\C}(\C)\; \equiv
    \; \C \smallsetminus \{0\}
\]
\end{corollary}

\begin{proof}
Every set contained into ${\rm Exp}_n(\C)$ (for any $n \in \N$)
of this form with bounded height (that is, $\R$-component) is
contained into a sub-helicoid ${\rm Exp}_n(B) \subseteq {\rm
Exp}_n(\C)$ of finite height, thus under its inversed image via
${\rm Exp}_n$ into a bounded strip $B$ of $\C$, for some $M > 0$.

The diffeomorphism comes from the previous \emph{Th.} \ref{The
Riemann surface} and the construction of injective complex
exponential on $\exp_{\C} (\C)$, as well as the remarks in the
beginning of the present \emph{Section} concerning the covering
space approach. We note that the covering space map coincides
with the immersion map into $\C \smallsetminus \{0\}$.
\end{proof}

\medskip
\noindent \textbf{(Part II) \underline{The Realization Theorem}.}

\medskip \medskip

\noindent We have not so far shown the explicit relation between
the surfaces ${\rm Exp}_a(\C)$ and $\Sigma_{log}$.

We recall that we manage to construct an injective $\exp_{\C}$
and thus a well-defined $\log_{\C}$ on ${\rm Exp}_a(\C) \subseteq
\R^3$ via convergence of helicoids. Taking into consideration the
initial \textit{"defining characterization" of B. Riemann himself}
on the surfaces of multi-valued holomorphic functions (as the one
of $\log_{\C}$), \emph{"surface onto which the function considered
becomes single-valued"} (see also \cite{Remm}), we conclude that
${\rm Exp}_a(\C)$ and $\Sigma_{log}$ must necessarily be
different representations of the "same" object, but the former
living into $\R^3$ instead of the latter which lives into $\R^4$.

This indeed turns to be the case, where "same" is expounded due
to the dimensional reduction to a realizable level (of
$3$-dimensions) as follows:

\begin{theorem} \textbf{(The Realization of $\Sigma_{log}$ into $\R^3$.)}
 \label{The Realization theorem}
The following diffeomorphism holds true
\[
    \Xi \;\;:\;\;\;\;{\rm Exp}_a (\C) \;\; {\xrightarrow{\;\;\;\;\;\; \cong \;\;\;\;\;\;}}
    \;\; \Sigma_{log}^{C^{\infty}}
\]
Thus, the (exponential) helicoid is $C^{\infty}$- diffeomorphic
to the Logarithmic Riemann surface, when the latter is equipped
with the induced smooth sub-atlas $\mathcal{A}_{\Sigma_{log}}
^{C^{\infty}}$ of the holomorphic $\mathcal{A}_{\Sigma_{log}}
^{\mathcal{O}}$ coming from $\R^4$.
\end{theorem}

The diffeomorphism $\Xi$ will be constructed in the proof. As a
consequence, \ref{The Realization theorem} implies that we can
equip the helicoid surface ${\rm Exp}_a (\C)$ with a holomorphic
structure via the bijective correspondence of it with
$\Sigma_{log}$, that is, if $(U, \phi)$ is a holomorphic chart of
$\mathcal{A}_{\Sigma_{log}}^{\mathcal{O}}$, then define an atlas
of ${\rm Exp}_a (\C)$ as follows
\[
  \mathcal{A}_{{\rm Exp}_a (\C)}^{\mathcal{O}} \;:=\; \{\bigl( \Xi^{-1}(U),
  \Xi^{-1} \circ \phi \bigr) \;/\;\; (U, \phi) \in \mathcal{A}_{\Sigma_{log}}
  ^{\mathcal{O}} \}
\]
Of course, this holomorphic atlas of the helicoid does
\underline{not} coincide with its induced smooth atlas coming
from the ambient $3$- space.
\begin{proof}
${\rm Exp}_a (\C)$ is given as the following subset of $\R^3$:
\[
  \{(e^u \cos v, e^u \sin v, a v)\;/\;\;\; u, \; v \in \R\}
\]
and $\Sigma_{log}$ as the subset of $\R^4$
\[
  \{(u, v, e^u \cos v, e^u \sin v)\;/\;\;\; u, \; v \in \R\}
\]
Define a map $\Xi$ as
\[
   \Xi \;:\;\; (e^u \cos v, e^u \sin v, a v) \; \longmapsto \;
    (u, v, e^u \cos v, e^u \sin v)
\]
Simple algebraic manipulations show that this map is a bijection
with inverse the projection onto $\R^3$ composed with a rigid
motion
\[
  \Omega \equiv \Xi^{-1} \;:\;\; (u, v, e^u \cos v, e^u \sin v)
   \; \longmapsto \;  (e^u \cos v, e^u \sin v, a v)
\]
both maps can be easily seen to have smooth components and this
completes the proof.
\end{proof}

\section{Multi-dimensional Generalizations.}

In this section we present the reasonable high-dimensional
analogues of the results presented in the previous
\emph{Sections}. Notwithstanding, the new context suggests that
the results will no longer stand in the region of Riemann
Surfaces, but will be of a more general differential-geometric
nature.

\noindent The following definitions arise naturally from the
respective \ref{Helicoid def}, \ref{Exp definition} of the
\emph{Section 1}.

\begin{definition} \textbf{(Multi-Exponential field)}
 \label{Multi-Exponential def}
The \emph{(exponential) multi-helicoid} is (in complex
coordinates) the map ($a_1, \ldots, a_m > 0$ parameters):
\[
   {\rm Exp}_{a_1, \ldots , a_m}\;:\;\;\;\;\; \C^m\; \larrow \;
   (\C \smallsetminus \{0\})^m \times \R^m \;\;\;:
\]
defined as
\[
   (z_1, \ldots , z_m) \; \longmapsto \;(\exp_{\C}(z_1), \ldots ,
    \exp_{\C}(z_m); \; a_1 Im (z_1), \ldots , a_m Im (z_m))
\]

that is, (up to a natural identitification $\C^m \times \R^m \cong
(\C \times \R)^m$)

\[
  {\rm Exp}_{a_1, \ldots , a_m}\;:=\; \prod_{1 \leq k \leq m} {\rm Exp}_{a_k}
\]

This is a vector field by \ref{Structural Properties of Exp}
(used component-wise), defined on the (flat) trivially embedded
submanifold $\C^m$ of $\R^{3m}$ and valued in (the tangent bundle
of) the ambient $\R^{3m}$ by \ref{Euclidean vector fields}.

Using \ref{Geometric Properties of Exp}, ${\rm Exp}_{a_1, \ldots
, a_m}(\C^m)$ is a globally coordinated smooth submanifold of
(real) dimension $\dim_{\R}( {\rm Exp}_{a_1, \ldots , a_m}(\C^m))
= 2m$. In fact, a \textbf{multi-helicoid}, in complete analogy
with \ref{Helicoid def}.
\end{definition}

The proofs of the following results go without saying using the
respective results of the previous sections, so we will not
bother showing anything explicitly. A component-wise argument
suffices in all cases, due to the product form of the definition
\ref{Multi-Exponential def}.

If we set $a_1 = \cdots = a_m = 1/ n$, we obtain a multi-sequence
of vector fields and in turn of multi-helicoid submanifolds of
$\R^{3m}$:
\[
    {\rm Exp}_{n, \ldots, n}(z_1, \ldots, z_m) \; := \; (\exp_{\C}(z_1), \ldots ,
    \exp_{\C}(z_m) ; \; \frac{1}{n} Im (z_1), \ldots ,  \frac{1}{n} Im (z_1))
\]

We are now in position to present the analogous result of Theorem
\ref{The Riemann surface}.

\begin{theorem} \textbf{(Convergence to Injective multi-Exponential)}
 \label{Convergence to injective multi Exp}

\noindent The sequence of exponential vector fields $({\rm
Exp}_{n, \ldots, n})_{n \in \; \N}$ converges
\underline{point-wise} to a $C^{\infty}$- diffeomorphism
$\underset{m}{\prod} \exp_{\C} \;:\;\; \C^m \; \larrow (\C
\smallsetminus \{0\})^m$ :
\[
\begin{picture}(300,220)
\put(240,165){\makebox{${\rm Exp}_{n-1, \ldots , n-1}(\C^m)$}}
  \put(25,165){\makebox{$\C^m$}}
  \put(42,170){\vector(1,0){195}}

\put(250,147){\makebox{$\scriptstyle{\cong}$}}
 \put(52,147){\makebox{$\scriptstyle{\cong}$}}

\put(211,125){\makebox{${\rm Exp}_{n, \ldots , n}(\C^m)$}}
  \put(60,125){\makebox{$\C^m$}}
  \put(77,130){\vector(1,0){131}}

 \put(234,121){\vector(-1,-1){18}}
\put(70,122){\vector(1,-1){20}} \put(93,95){\makebox{$\cdot$}}
\put(98,90){\makebox{$\cdot$}} \put(103,85){\makebox{$\cdot$}}

\put(210,95){\makebox{$\cdot$}} \put(205,90){\makebox{$\cdot$}}
\put(200,85){\makebox{$\cdot$}}

 \put(280,160){\vector(0,-1){100}}
  \put(30,160){\vector(0,-1){100}}
\put(280,46){\vector(0,-1){35}}
\put(270,30){\makebox{$\scriptstyle{\cong}$}}

 \put(271,160){\vector(-1,-1){25}}
\put(35,160){\vector(1,-1){25}}

\put(-8,200){\makebox{$\cdot$}} \put(313,200){\makebox{$\cdot$}}
\put(-13,205){\makebox{$\cdot$}} \put(318,205){\makebox{$\cdot$}}
\put(-18,210){\makebox{$\cdot$}} \put(323,210){\makebox{$\cdot$}}
 \put(-2,198){\vector(1,-1){25}}
 \put(307,198){\vector(-1,-1){25}}

 \put(60,121){\vector(-1,-4){15}}
 \put(241,121){\vector(1,-4){15}}

\put(22,110){\makebox{$\scriptstyle{\Phi}$}}
 \put(45,90){\makebox{$\scriptstyle{\Phi}$}}
\put(33,110){\makebox{$\scriptstyle{\cong}$}}
 \put(57,90){\makebox{$\scriptstyle{\cong}$}}
 \put(240,90){\makebox{$\scriptstyle{\Psi}$}}
\put(285,110){\makebox{$\scriptstyle{\Psi}$}}

 \put(183,0){\makebox{$\underset{m}{\prod} \exp_{\C}(\C) = (\C \smallsetminus \{0\})^m$}}

\put(75,53){\vector(1,0){165}}
 \put(7,49){\makebox{$\underset{\larrow}{\lim}\; \C^m \cong \C^m$}}
\put(242,49){\makebox{$\underset{\larrow}{\lim} \; {\rm Exp}_{n,
\ldots, n}(\C^m)$}}
   \put(38,45){\vector(4,-1){142}}
\put(136,60){\makebox{$\scriptstyle{\underset{\larrow}{\lim}\;
{\rm Exp}_{n, \ldots, n}}$}}
 \put(145,45){\makebox{$\scriptstyle{\cong}$}}

\put(130,134){\makebox{$\scriptstyle{{\rm Exp}_{n, \ldots, n}}$}}
\put(124,174){\makebox{$\scriptstyle{{\rm Exp}_{n-1, \ldots,
n-1}}$}} \put(135,122){\makebox{$\scriptstyle{\cong}$}}
\put(138,162){\makebox{$\scriptstyle{\cong}$}}
\put(104,17){\makebox{$\scriptstyle{\underset{m}{\prod}
\exp_{\C}}$}}
\end{picture}
\]
\noindent The convergence is \underline{uniform} on the bounded
multi-strips $B^m \equiv \{z_k \in \C \;:\; |Im (z_k)|< M_k
\;/\;\; M_k
> 0,\; 1 \leq k \leq m\} \; \subseteq \; \C^{m}$.
\[
    {\rm Exp}_{n, \ldots, n}|_{B^m} \; \underset{U}
    {\xrightarrow{\;\;\;\;n \larrow \infty\;\;\;\;}}\;
    \underset{m}{\prod} \exp_{\C} \times \{0\}|_{B} \;
     \equiv \; \underset{m}{\prod} \exp_{\C}|_{B}
\]
\end{theorem}

\medskip

\begin{remark}
The methods we expounded for these specific construction raise the
question if and how they may be extended and applied to other
surfaces, in the framework of a general approach that would, at
least, be in position to give back the already known facts
concerning the classical Riemann surfaces, e.g. those of
$\log_{\C}(z)$ and $\sqrt[n]{z}$ which, the latter, more or less
constituted the beginnings of the subject of the study of
holomorphic $1$-manifolds.

\noindent Our problem is focused in the quest of an appropriate
sequence of surfaces in $\R^3$ (if any), whereon the complex
function in question (properly extended as a vector field) will
become single valued, and the covering surface will be obtained in
a uniform limit.

The above attractive concept might be of significant importance in
the finite-sheeted Riemann surfaces of algebraic functions
satisfying the general formula $a_n(z) [f(z)]^n +
a_{n-1}(z)[f(z)]^{n-1}+ \cdots + a_0(z) = 0$.

\medskip

\textit{In a nutshell, recovering via differential geometry pure
analytic information concerning complex functions through their
Riemann surfaces would give a measure of the amount that (the
latter) in fact depend on the holomorphic structure of the plane
and, more generally, of every complex manifold modeled upon it}.
\end{remark}

\medskip

\noindent {\bf Acknowledgement.} The author is profoundly indebted
to prof. \textit{E. Vassiliou}, \textit{A. Kartsaklis} \&
\textit{P. Krikelis} of University of Athens, Section of Algebra
and Geometry for their multifarious contribution in many levels
during the preparation of this paper.


\begin{thebibliography}{16}
\bibitem{B-N} J. Bak, D. J. Newman, \emph{Complex Analysis},
Spriger-Verlag, New York, 1982.

\bibitem{Behnke-Stein} H. Behnke and K. Stein, \emph{Entwicklung
 analytischer Funktionen auf Riemannschen Fl\"achen},
Math. Ann., 120 (1947-49), pp. 430  461.

\bibitem{Bre} G. E. Bredon, \emph{Topology and Geometry},
Springer, 1997.

\bibitem{B-C} Brickell F., Clark R.S., \emph{Differential
Manifolds}, Van Nostrand, Reinhold Company.

\bibitem{DC} Do Carmo, M.P., \emph{Riemannian Geometry},
Birkha\"{u}ser, 1992.

\bibitem{Eis} L. Pf. Eisenhart, \emph{Riemannian Geometry},
Princeton University Press, 1964.

\bibitem{G-H-L} S. Gallot, D. Hulin, J. Lafontaine,
\emph{Riemannian Geometry}, Springer-Verlag, 1993, 2nd printing.

\bibitem{Gunning-Nar} R. C. Gunning and R. Narasimhan \emph{Immersion
 of open Riemann surfaces}, Math. Ann., 174 (1967), pp.
103  108.

\bibitem{F-K} Farkas H.M., Kra I., \emph{Riemann Surfaces}, New
York, Springer-Verlag, 1980.

\bibitem{Hormander} L. H\"{o}rmander, \emph{An Introduction to Complex Analysis
 in Several Varables}, North-Holland Publishing, 1979.

\bibitem{K-N} Kobayashi S., Nomizu K., \emph{Foundations of
Differential Geometry I, II},1963, 1969, Interscience.

\bibitem{Krasnov} K. Krasnov, \emph{Black-hole thermodynamics and Riemann surfaces},
 Class. Quantum Grav. 20 2235 - 2250, 2003.


\bibitem{L-S-Y} K. liu, X. Sun, S.-T. Yau, \emph{Geometric aspects of the Moduli space of Riemann surfaces},
 arXiv: math.DG/0411247, pre-print.

\bibitem{Mac Lane} S. Mac Lane, \emph{Categories for the Working
Mathematician}, Springer - Verlag, 1971.

\bibitem{Mun} J. R. Munkres, \emph{Topology, a first course},
Prentice-Hall, Inc., Englewood Cliffs, New Jersey.

\bibitem{Nar} R. Narasimhan, \emph{Analysis on Real and Complex
Manifolds}, Masson and Cie, Paris, North Holland, 1968.

\bibitem{ON} O' Neil, \emph{Elementary Differential Geometry}, New
York and London, Academic Press, 1997.

\bibitem{Remm} R. Remmert, \emph{From Riemann Surfaces to Complex
Spaces}, Lecture notes, (Nice), french mathematical society, 1998.

\bibitem{Sharipov}, R. A. Sharipov, \emph{Course of Differential Geometry, the
textbook}, Russian Federal Committee for Higher Education, Bashkir
State University, 1996.

\bibitem{S} Saunders D.J., \emph{The Geometry of Jet Bundles},
Cambridge University Press, 1989.

\bibitem{S-T} Singer I., Thorpe J.A., \emph{Lecture notes on
Emementary Topology and Geometry}, Scott, Foresman and Co, 1967.

\bibitem{Spr} Springer G.,\emph{Introduction to Riemann Surfaces},
Chelsea Publishing Company, 1981.

\end{thebibliography}
\end{document}